\documentclass[12pt]{amsart}

\title{M{\"o}bius disjointness for $C^{1 + \eps}$ skew products} 

\author{Alexandre de Faveri}
\address{Department of Mathematics, Caltech, 1200 E. California Blvd., Pasadena, CA 91125, USA}
\email{\href{mailto:afaveri@caltech.edu}{afaveri@caltech.edu}}

\usepackage{amssymb, amsmath, amsthm}
\usepackage{hyperref, enumitem}
\usepackage[space]{cite}
\usepackage[centering]{geometry}
\usepackage{xcolor, graphicx}

\newcommand{\dd}{\, d}
\newcommand{\R}{\mathbb{R}}

\newcommand{\Z}{\mathbb{Z}}
\newcommand{\Q}{\mathbb{Q}}
\newcommand{\T}{\mathbb{T}}

\newcommand{\eps}{\varepsilon}

\newtheorem{theorem}{Theorem}
\newtheorem{lemma}{Lemma}

\newtheorem{proposition}{Proposition}

\newtheorem{corollary}{Corollary}

\newtheorem{remark}{Remark}

\newenvironment{boldproof}[1][\proofname]{\proof[\bf #1:]}{\endproof} 

\DeclareMathOperator{\Var}{Var}
\DeclareMathOperator{\Id}{Id}

\begin{document}

\begin{abstract}
	We show that for $\varepsilon > 0$, every $C^{1+\varepsilon }$ skew product on $\mathbb{T}^2$ over a \mbox{rotation} of $\mathbb{T}^1$ satisfies Sarnak's conjecture. This is an improvement of earlier results of Ku{\l}aga-Przymus-Lema\'{n}czyk, Huang-Wang-Ye, and Kanigowski-Lema{\'n}czyk-Radziwi{\l}{\l}.
\end{abstract}

\maketitle


\section{Introduction}

Let $(X, d)$ be a compact metric space and $T:X \to X$ be a homeomorphism. If the topological dynamical system $(X, T)$ has (topological) entropy zero, then Sarnak's conjecture \cite{sarnak_conj, sarnak_conj_online} predicts that

\begin{equation*}\label{disjointness_eq}
    \lim_{N \to \infty} \frac{1}{N} \sum_{n \leq N} f(T^n x) \mu(n) = 0
\end{equation*}
for any continuous $f:X \to \R$ and every $x \in X$. When this holds, we say that the system $(X, T)$ is \textit{M{\"o}bius disjoint}. 

Sarnak's conjecture has been proved for a variety of dynamical systems: see for instance \cite{mand_rivat, rank_one, bounded_depth, davenport, fran_host, rank_one_pow, green_tao, prime_digits, digital_functions, homogeneous, BSZ}. A common feature of all the results listed is that the underlying system is \textit{regular}, in the sense that for every $x \in X$ the sequence $\frac{1}{N} \sum_{n \leq N} \delta_{T^n(x)}$ converges in the weak-* topology to some $T$-invariant Borel probability measure on $X$.

Let $\T := \R/\Z$ denote the circle. In this paper we will deal with the so-called \textit{Anzai skew products} $(\T^2, T_{\alpha, \phi})$, where $\alpha \in \R$, $\phi:\T \to \T$ is a continuous map and the transformation is given by
\begin{equation*}
    T_{\alpha, \phi}(x, y) := (x + \alpha, y + \phi(x))
\end{equation*}
for all $(x, y) \in \T^2$. We often denote the system simply by $T_{\alpha, \phi}$.

Observe that $T_{\alpha, \phi}$ is \textit{distal}, so it has zero topological entropy and therefore we expect it to be M\"obius disjoint. In fact, these skew products are the basic building blocks in Furstenberg's classification of minimal distal flows \cite{furstenberg_distal_flows}, so understanding them is the first step towards establishing Sarnak's conjecture for this important general case. The main novel dynamical challenge that arises when one deals with skew products is that they provide some of the simplest examples of irregular dynamics. Indeed, Furstenberg \cite{furstenberg_irregular} showed that $T_{\alpha, \phi}$ is not regular for some $\alpha$ and some analytic $\phi$.

Lifting $\phi:\T \to \T$ to the real line, we can write $\phi(x) = cx + \Tilde{\phi}(x)$ for all $x \in \T$, where $c \in \Z$ is the \textit{topological degree} of $\phi$ and $\Tilde{\phi}:\T \to \R$ is a continuous $1$-periodic function, unique up to shifts by $\Z$ (we fix an arbitrary choice). Ku{\l}aga-Przymus and Lema\'{n}czyk \cite{C1_aa} have shown that if $\phi \in C^{1+\eps}$ for some $\eps > 0$ then $T_{\alpha, \phi}$ is M\"obius disjoint for a topologically generic set of $\alpha$. Furthermore, they proved M\"obius disjointness of $T_{\alpha, \phi}$ when $\alpha \in \Q$, assuming only continuity of $\phi$ \cite[Proposition 2.3.1]{C1_aa}, so from now on we assume $\alpha \in \R \setminus \Q$. A further consequence of their work \cite[Remark 2.5.7]{C1_aa} (see also \cite[Corollary 2.6]{wang_analytic}) is that if $\phi$ is assumed to be Lipschitz continuous, then Sarnak's conjecture holds for $T_{\alpha, \phi}$ whenever $c \not= 0$. Therefore, with the underlying assumption on $\phi$ in mind, we can deal only with topological degree zero from now on, and with an abuse of notation we identify $\phi$ with $\Tilde{\phi}$. 

The first M\"obius disjointness result for all $\alpha$ was established by Liu and Sarnak \cite{liu_sarnak}, who proved it for $\phi$ analytic and satisfying the technical condition $\widehat{\phi}(m) \gg e^{-\tau|m|}$ for some $\tau > 0$. This was the first time Sarnak's conjecture was proved for a system that is not regular (since Furstenberg's example satisfies the condition). A refinement of this result was recently obtained by Wang \cite{wang_analytic}, who removed the need for a lower bound on Fourier coefficients, obtaining M\"obius disjointness of $T_{\alpha, \phi}$ for all analytic $\phi$. Huang, Wang, and Ye \cite{C_infty} later improved this to cover all $\phi \in C^\infty$. Finally, using the work of Matom\"{a}ki and Radziwi{\l}{\l} \cite{M-R} on the behavior of $\mu$ in short intervals, Kanigowski, Lema{\'n}czyk, and Radziwi{\l}{\l} \cite{C2} established M\"obius disjointness of $T_{\alpha, \phi}$ for all $\phi \in C^{2 + \eps}$ subject to the condition $\widehat{\phi}(0) = 0$, where $\eps > 0$ is arbitrary.

Our main result is a simultaneous improvement of the works of Ku{\l}aga-Przymus-Lema\'{n}czyk \cite{C1_aa}, Huang-Wang-Ye \cite{C_infty}, and Kanigowski-Lema{\'n}czyk-Radziwi{\l}{\l} \cite{C2}:

\begin{theorem}\label{main_thm}
    Let $\eps > 0$. For any $\alpha \in \R$ and $\phi:\T \to \T$ of class $C^{1 + \eps}$, the skew product $T_{\alpha, \phi}$ is M{\"o}bius disjoint.
\end{theorem}

The proof follows the ideas laid out by Kanigowski-Lema{\'n}czyk-Radziwi{\l}{\l} in \cite{C2}, but instead of aiming for a polynomial rate of convergence for $T_{\alpha, \phi}^{r_n} \to \Id$ in the uniform norm (along some unbounded sequence $\{r_n\}_{n \geq 1}$), we establish a polynomial rate of convergence for $T_{\alpha, \phi}^{r_n} \to \Id$ in the $L^2(\nu)$ norm, for each $T_{\alpha, \phi}$-invariant Borel probability measure $\nu$. The difficulties in dealing with every such $\nu$ are overcome because they all project to the Lebesgue measure in the first coordinate. We also remove the condition $\widehat{\phi}(0) = 0$ present in \cite{C2} by slightly modifying their choice of the sequence $\{r_n\}_{n\geq 1}$. 

Another important ingredient is better control of some sums related to the Fourier coefficients of $\phi$, where the Diophantine properties of $\alpha$ play an important role. The idea here is that not many $q$'s at a given scale can make $\Vert q \alpha \Vert$ small (i.e. be denominators of good rational approximations of $\alpha$). Furthermore, the $q$'s at a given scale that give rise to rational approximations of similar quality must be somewhat well-spaced. We apply the Denjoy-Koksma inequality to appropriately chosen functions in order to extract that information (see \autoref{cont_frac_section}).

The smoothness exponent $1 + \eps$ seems to be the limit of this argument. Indeed, we prove in \autoref{counterexample_section} that if one only assumes that $\phi \in C^1$ then, at least along the sequence of best rational approximations of the irrational $\alpha$, the rate of rigidity of $T_{\alpha, \phi}$ can be logarithmic even when $\widehat{\phi}(0) = 0$.

In \autoref{abs_con_section} we show that our ideas can be used to extend some general rigidity results so far only known for functions of mean zero to the general case. A modification of \autoref{PR_lemma} to obtain uniform polynomial rates of rigidity in the case $\phi \in C^{1+\eps}$ is also discussed.

Finally, in \autoref{flows_rokhlin_section} we use our argument to deduce new M{\"o}bius disjointness results for flows in $\T^2$ and Rokhlin extensions.

\subsection*{Notations}
    For a topological dynamical system $(X, T)$, let $M(X, T)$ be the set of $T$-invariant Borel probability measures on $X$. Write $\Vert \cdot \Vert$ for the distance to the nearest integer (which we use as the metric in $\T$), $d(\cdot, \cdot)$ for the product metric in $\T \times \T$, corresponding to $\Vert \cdot \Vert$ in each coordinate, and $\Vert \cdot \Vert_{L^2(\nu)}$ for the usual $L^2$ norm with respect to a measure $\nu$. We also abbreviate $e(x) := e^{2 \pi i x}$ and use the asymptotic notation $f(x) \ll g(x)$ (respectively $f(x) \ll_p g(x)$) to mean that there exists $C>0$ absolute (respectively depending only on the parameter $p$) such that $|f(x)| \leq C |g(x)|$ for all $x$ in the relevant range. Furthermore, $f(x) \asymp g(x)$ means $f(x) \ll g(x) \ll f(x)$.
     
\subsection*{Acknowledgments}
    I would like to thank my PhD advisor, Maksym Radziwi{\l}{\l}, for introducing me to this problem and for general advice and encouragement. Thanks also to Adam Kanigowski and Mariusz Lema{\'n}czyk for pointing out a nice simplification to my initial proofs of Lemmas \ref{cont_frac_whole_lemma} and \ref{cont_frac_slice_lemma}, and for providing valuable comments and references. I'm grateful to the American Institute of Mathematics (AIM) for their 2018 workshop on ``Sarnak's Conjecture'', which played a role in motivating this work.


\section{Reduction of \autoref{main_thm} to a rigidity result}

As previously outlined, we can assume that $\alpha \in \R \setminus \Q$ and $\deg(\phi) = 0$, so $\phi$ can be realized as a function from $\T$ to $\R$ of class $C^{1 + \eps}$, which by an abuse of notation we still denote by $\phi$. Observe that $\phi$ is in particular Lipschitz continuous, so we have pointwise convergence of its Fourier series, and the smoothness condition gives

\begin{equation}\label{fourier_eq}
    \phi(x) = \sum_{q \in \Z} c_q e(qx) \quad \text{with}\quad c_q \ll_\phi \frac{1}{|q|^{1 + \eps}} \ \text{for} \  q \not= 0.
\end{equation}

The key to the proof of \autoref{main_thm} is the result below, which is motivated by \cite{C2}.

\begin{lemma}\label{PR_lemma}
    Let $0 < \eps < \frac{1}{100}$ and $\alpha \in \R \setminus \Q$. If $\phi: \T \to \R$ is of class $C^{1 + \eps}$, then there exists an unbounded sequence of positive integers $\{r_n\}_{n\geq1}$ such that
    \begin{equation*}
        \int_{\T \times \T} d(T_{\alpha, \phi}^{r_n}(x, y), (x, y))^2 \dd\nu(x, y) \ll_\phi r_n^{-\eps/100}
    \end{equation*}
    for any $\nu \in M(\T^2, T_{\alpha, \phi})$, where the implied constant does not depend on $\nu$.
\end{lemma}

Assuming \autoref{PR_lemma}, we can easily prove \autoref{main_thm}:

\begin{boldproof}[Proof of \autoref{main_thm}] Let $\{r_n\}_{n\geq 1}$ be the sequence from \autoref{PR_lemma}. For any $\nu \in M(\T^2, T_{\alpha, \phi})$, continuous $f:\T^2 \to \R$, and $k\in \Z$, the triangle inequality and the $T_{\alpha, \phi}$-invariance of $\nu$ imply 
\begin{equation}\label{triangle_ineq_eq}
    \Vert f \circ T_{\alpha, \phi}^{k r_n} - f \Vert_{L^2(\nu)}^2 \leq |k| \sum_{j = 1}^{|k|} \Vert f \circ T_{\alpha, \phi}^{j r_n} - f \circ T_{\alpha, \phi}^{(j-1) r_n} \Vert_{L^2(\nu)}^2 = k^2 \cdot \Vert f \circ T_{\alpha, \phi}^{r_n} - f \Vert_{L^2(\nu)}^2.
\end{equation}

If $f$ is also Lipschitz continuous, then using \autoref{PR_lemma} we get
\begin{equation}\label{pointwise_norm_eq}
    \Vert f \circ T_{\alpha, \phi}^{r_n} - f \Vert_{L^2(\nu)}^2 \ll_f \int_{\T \times \T} d(T_{\alpha, \phi}^{r_n}(x, y), (x, y))^2 \dd\nu(x, y) \ll_\phi r_n^{-\eps/100}.
\end{equation}

Therefore, \eqref{triangle_ineq_eq} and \eqref{pointwise_norm_eq} together give
\begin{equation*}
    \lim_{n\to \infty} \sum_{|k| \leq r_n^{\eps/400}} \Vert f \circ T_{\alpha, \phi}^{k r_n} - f\Vert_{L^2(\nu)}^2 = 0
\end{equation*}
    for every $\nu \in M(\T^2, T_{\alpha, \phi})$, which is precisely the \textit{PR rigidity} condition of \cite{C2} (using the linearly dense family $\mathcal{F}$ of Lipschitz continuous functions) for the system $(\T^2, T_{\alpha, \phi})$, so \cite[Theorem 1.1]{C2} implies M{\"o}bius disjointness for this skew product, and \autoref{main_thm} is proved.
    
\end{boldproof}


\section{Continued fractions and some arithmetic estimates}\label{cont_frac_section}

Before proceeding to the proof of \autoref{PR_lemma}, we recall some properties of continued fractions. Let $\frac{p_n}{q_n}$, with $q_n > 0$ and $(p_n, q_n) = 1$, be the $n$-th convergent of the continued fraction expansion $[a_0;a_1, a_2, \dots]$ of the irrational $\alpha$, so that $a_i \geq 1$ for $i \not= 0$. Then:

\begin{enumerate}[label={(P\arabic*)}]
    \item \label{P1_item} $q_0 = 1$, $q_1 = a_1$ and $q_{n+1} = a_{n+1} q_n + q_{n-1}$ for $n \geq 1$;\label{cde}
    \item \label{P2_item} $\frac{1}{q_{n+1} + q_n} < \Vert q_n \alpha \Vert < \frac{1}{q_{n+1}}$;
    \item \label{P3_item} If $0 < q < q_{n+1}$, then $\Vert q_n \alpha \Vert \leq \Vert q \alpha \Vert$.
\end{enumerate}

The main technical tool that allows us to quickly explore the Diophantine properties of $\alpha$ through its continued fraction is the following inequality:

\begin{proposition}[Denjoy-Koksma inequality]\label{denjoy_koksma_prop} Let $\alpha \in \R \setminus \Q$. If $f:\T \to \R$ is of bounded variation, which we denote by $\Var(f)$, then for any $n \geq 0$ and $x \in \T$ we have
\begin{equation*}
    \left|\sum_{j = 0}^{q_n-1} f(x + j \alpha) - q_n \int_{\T} f(z) \dd z \right| \leq \Var(f).
\end{equation*}
\end{proposition}

\begin{boldproof}[Proof of \autoref{denjoy_koksma_prop}]
    See \cite[VI.3.1]{herman_abs_con}.
    
\end{boldproof}

The next two lemmas encapsulate estimates related to continued fractions that will be necessary to prove \autoref{PR_lemma}.

\begin{lemma}\label{cont_frac_whole_lemma}
    For any $\alpha \in \R \setminus \Q$ and $k \geq 2$,
    \begin{equation*}
        \sum_{0 < |q| < q_k} \frac{1}{\Vert q \alpha \Vert^{2}} \asymp q_k^2.
    \end{equation*}
\end{lemma}

\begin{boldproof}[Proof of \autoref{cont_frac_whole_lemma}]
    The lower bound comes from positivity and the single term $q = q_{k-1}$, by \ref{P2_item}. The upper bound follows from \cite[Lemma 2.5]{koksma_banach} (see also \cite[Lemma 1]{rigidity_L_2} for a partial result). We give a quick proof for completeness.
    
    Assume $0 < q < q_k$, as the sum over negative $q$ is the same. Consider $f:\T \to \R$ given by
    \begin{equation*}
        f(z) = 
        \begin{cases}
        (2q_k)^2, & \text{if} \ \Vert z \Vert \leq \frac{1}{2q_k} \\
        \frac{1}{\Vert z \Vert^2}, & \text{if} \ \Vert z \Vert > \frac{1}{2q_k}
        \end{cases}.
    \end{equation*}
    Observe that by \ref{P2_item} and \ref{P3_item}, $\Vert q \alpha \Vert > \frac{1}{2 q_k}$ for all $0 < q < q_k$, so by the Denjoy-Koksma inequality we conclude that
    \begin{equation*}
    \begin{split}
        \sum_{0 < q < q_k} \frac{1}{\Vert q \alpha \Vert^{2}} = \sum_{q = 1}^{q_k - 1} f(q \alpha) &\leq |f(0)| + q_k \left| \int_\T f(z) \dd z \right| + \Var(f) \\
        & = 4 q_k^2 + q_k(8q_k - 4) + (8 q_k^2 - 8) \ll q_k^2,
    \end{split}
    \end{equation*}
    as desired.
    
\end{boldproof}

\begin{lemma}\label{cont_frac_slice_lemma}
    For any $\alpha \in \R \setminus \Q$, $k \geq 1$, and $1 \leq c \leq q_k$, 
    \begin{equation*}
        \sum_{q_k \leq |q| < q_{k+1}} \frac{1}{q^2} \min\left\{\frac{1}{\Vert q \alpha \Vert^{2}}, c^2\right\} \ll \frac{c}{q_k}.
    \end{equation*}
\end{lemma}

\begin{boldproof}[Proof of \autoref{cont_frac_slice_lemma}] 
    We can assume $q_{k} \leq q < q_{k+1}$ since the sum over negative $q$ is the same. 
    
    Consider $f:\T \to \R$ given by
    \begin{equation*}
        f(z) = 
        \begin{cases}
        c^2, & \text{if} \ \Vert z \Vert \leq \frac{1}{c} \\
        \frac{1}{\Vert z \Vert^2}, & \text{if} \ \Vert z \Vert > \frac{1}{c}
        \end{cases}.
    \end{equation*}
    Observe that $f(q\alpha) = \min\left\{\frac{1}{\Vert q \alpha \Vert^{2}}, c^2\right\}$, so \ref{P1_item} gives
    \begin{equation*}
    \begin{split}
        &\sum_{q_k \leq q < q_{k+1}} \frac{1}{q^2} \min\left\{\frac{1}{\Vert q \alpha \Vert^{2}}, c^2\right\} = \sum_{j = 1}^{a_{k+1} - 1}\sum_{j q_k \leq q < (j+1) q_{k}} \frac{f(q\alpha)}{q^2} + \sum_{a_{k+1} q_k \leq q < q_{k+1}} \frac{f(q\alpha)}{q^2} \\
        & \leq \sum_{j = 1}^{a_{k+1} - 1}\frac{1}{(jq_k)^2} \sum_{r=0}^{q_{k}-1} f(j q_{k} \alpha + r\alpha) + \frac{4}{q_{k+1}^2} \sum_{r=0}^{q_{k-1}-1} f(a_{k+1}q_k \alpha + r\alpha).
    \end{split}
    \end{equation*}
    
    Using the Denjoy-Koksma inequality for the sums over $r$ as in the proof of \autoref{cont_frac_whole_lemma}, since $\int_\T f(z) \dd z \asymp c$ and $\Var(f) \ll c^2$ by direct computation, we conclude that the remaining expression is
    \begin{equation*}
        \ll \sum_{j=1}^\infty \frac{q_k c + c^2}{(jq_k)^2} + \frac{q_{k-1}c + c^2}{q_{k+1}^2} \ll \frac{c}{q_k},
    \end{equation*}
    so we are done.
    
\end{boldproof}


\section{Polynomial rate of rigidity in \texorpdfstring{$C^{1+ \eps}$}{C 1+e}: the proof of \autoref{PR_lemma}}

At last, we are ready to prove our main lemma.

\begin{boldproof}[Proof of \autoref{PR_lemma}]
    Denoting $S_r(g)(x) := g(x) + g(x+\alpha) + \dots + g(x + (r-1)\alpha)$, we have
    \begin{equation*}
        T_{\alpha, \phi}^{r_n}(x, y) = (x + r_n \alpha, y + S_{r_n}(\phi)(x)),
    \end{equation*}
    so that
    \begin{equation*}
        d(T_{\alpha, \phi}^{r_n}(x, y), (x, y))^2 \asymp \Vert r_n \alpha \Vert^2 + \Vert S_{r_n}(\phi)(x) \Vert^2.
    \end{equation*}
    Therefore,
    \begin{equation}\label{intermed_eq}
    	D_n := \int_{\T \times \T} d(T_{\alpha, \phi}^{r_n}(x, y), (x, y))^2 \dd\nu(x, y) \asymp \Vert r_n \alpha \Vert^2 + \int_{\T \times \T} \Vert S_{r_n}(\phi)(x)\Vert^2 \dd \nu(x, y).	
    \end{equation}
    
    Consider the projection map $\pi(x, y) = x$. Observe that the integrand in \eqref{intermed_eq} is independent of the second coordinate, so we can rewrite the integral as
    \begin{equation}\label{integral_eq}
    \begin{split}
    	\int_{\T \times \T} \Vert S_{r_n}(\phi) (\pi (x, y))\Vert^2 \dd \nu(x, y) = \int_{\T} \Vert S_{r_n}(\phi)(x)\Vert^2 \dd (\pi_* \nu)(x).
    \end{split} 
    \end{equation}
    Since $\pi:(\T^2, T_{\alpha, \phi}) \to (\T, R_\alpha)$ is a map of topological dynamical systems (where in the image the transformation is $R_\alpha(x) := x + \alpha$) and $\nu$ is $T_{\alpha, \phi}$-invariant, the Borel probability measure $\pi_*\nu$ is $R_\alpha$-invariant. But $\alpha$ is irrational, so $(\T, R_\alpha)$ is uniquely ergodic and we conclude that $\pi_* \nu$ is the Lebesgue measure on $\T$.

Using the Fourier expansion of $\phi$ we get $S_{r_n}(\phi)(x) = \sum_{q \in \Z} c_q S_{r_n}(e_q)(x)$, where $e_q(x) := e(q x)$. A computation shows that 
\begin{equation*}
    S_{r_n}(e_q)(x) = e(qx)\frac{1 - e(q r_n \alpha)}{1 - e(q \alpha)}
\end{equation*} 
for $q \not= 0$ and $S_{r_n}(e_0)(x) = r_n$, so we can plug this into \eqref{integral_eq} and conclude, using the triangle inequality and replacing $\Vert \cdot \Vert$ by absolute values, that the integral there is bounded by a constant multiple of
\begin{equation}\label{square_eq}
	 \Vert c_0 r_n \Vert^2 + \int_\T \left| \sum_{q \not= 0} c_q \frac{1 - e(q r_n \alpha)}{1 - e(q \alpha)} e(qx) \right|^2 dx = \Vert c_0 r_n \Vert^2 + \sum_{q \not= 0} |c_q|^2 \left| \frac{1 - e(q r_n \alpha)}{1 - e(q \alpha)} \right|^2,
\end{equation}
where we have used Parseval for $S_{r_n}(\phi) - c_0 r_n \in L^2(\T)$.

Now, we make a preliminary choice of the sequence $\{r_n\}_{n\geq 1}$ by letting $r_n := \ell_n q_n$, where $q_n$ is the denominator of the $n$-th convergent of the continued fraction expansion of $\alpha$, as before, and $\ell_n \in \Z$ is chosen so that 
\begin{equation*}
    0 < \ell_n \leq q_n^\delta \quad \text{and} \quad \Vert \ell_n q_n c_0 \Vert < q_n^{- \delta},
\end{equation*} 
where $\delta := \eps/10$. Such $\ell_n$ exist for all $n$, by the Dirichlet approximation theorem.

Let $\lambda := \eps/100$. In what follows it is worth keeping in mind the rough hierarchy ``$\lambda \lll \delta \lll \eps$" behind our choice of parameters. We wish to show that $D_n \ll_\phi r_n^{-\lambda}$. With our choice of $\{r_n\}_{n\geq 1}$ the first term in the RHS of \eqref{intermed_eq} contributes at most
\begin{equation}\label{term_1_eq}
	\ell_n^2 \cdot \Vert q_n \alpha \Vert^2 < q_n^{2\delta} q_{n+1}^{-2} < q_n^{2\delta - 2} < q_n^{-\lambda(1+\delta)} \leq r_n^{-\lambda}, 
\end{equation}
so it is harmless. The first term of \eqref{square_eq} contributes
\begin{equation}\label{term_2_eq}
	\Vert c_0 \ell_n q_n \Vert^2 < q_n^{-2\delta} < q_n^{-\lambda(1+\delta)} \leq r_n^{-\lambda},
\end{equation}
and it is also harmless. 

We break the remaining terms into two parts, corresponding to $0 < |q| < q_n$ and $|q| \geq q_n$. Observe that $|1 - e(q\alpha)| \asymp \Vert q \alpha \Vert$ and $|1 - e(q r_n \alpha)| \leq 2$. Furthermore, $|S_{r_n}(e_q)(x)| \leq r_n$ by a trivial bound, so
\begin{equation}\label{term_3_eq}
\begin{split}
	\sum_{|q| \geq q_n} |c_q|^2 \left| \frac{1 - e(q r_n \alpha)}{1 - e(q \alpha)} \right|^2 &\ll_\phi \sum_{|q| \geq q_n} \frac{1}{|q|^{2+2\eps}} \min\left\{\frac{1}{\Vert q \alpha \Vert^{2}}, r_n^2\right\} \\
	&< q_n^{-2\eps} \ell_n^2 \sum_{k=n}^\infty \sum_{q_k \leq |q| < q_{k+1}} \frac{1}{q^{2}} \min\left\{\frac{1}{\Vert q \alpha \Vert^{2}}, q_n^2\right\} \\ 
	&\ll q_n^{-2\eps + 2\delta}\sum_{k = n}^\infty \frac{q_n}{q_k} \ll q_n^{-2\eps + 2\delta} < q_n^{-\lambda(1+\delta)} \leq r_n^{-\lambda},
\end{split}
\end{equation}
where we have used \eqref{fourier_eq}, \autoref{cont_frac_slice_lemma} (for $c = q_n \leq q_k$) and the fact that $q_{k+2} > 2q_k$ by \ref{P1_item}, so the $q_k$ grow exponentially.

It remains to deal with $0 < |q| < q_n$. In this case, we use $|1 - e(q\alpha)| \asymp \Vert q \alpha \Vert$ and $|1 - e(q r_n \alpha)|^2 \asymp \Vert q \ell_n q_n \alpha\Vert^2 \leq q^2 \ell_n^2 \cdot \Vert q_n \alpha \Vert^2 < q^2 q_n^{2\delta} q_{n+1}^{-2}$, so those terms contribute
\begin{equation}\label{term_4_eq}
\begin{split}
	\sum_{0 < |q| < q_n} |c_q|^2 \left| \frac{1 - e(q r_n \alpha)}{1 - e(q \alpha)} \right|^2 &\ll_\phi \frac{q_n^{2\delta}}{q_{n+1}^{2}} \sum_{0 < |q| <q_n} \frac{1}{|q|^{2\eps}} \frac{1}{\Vert q \alpha \Vert^{2}}.
\end{split}
\end{equation}

To deal with the sum over $q$ we consider two cases:

\begin{itemize}
	\item[\textbf{\textit{Case 1:}}] \textit{There is a subsequence $\{q_{b_n}\}_{n \geq 1}$ of $\{q_n\}_{n \geq 1}$ such that $q_{b_n + 1} \geq q_{b_n}^2$ for all $n \geq 1$.}
	
	In this case we take the subsequence $\{r_{b_n}\}_{n \geq 1}$ instead of the original sequence $\{r_n\}_{n \geq 1}$. Observe that \eqref{term_1_eq}, \eqref{term_2_eq} and \eqref{term_3_eq} still hold along any subsequence. In \eqref{term_4_eq} we can use the given condition and \autoref{cont_frac_whole_lemma} to get the upper bound 
	\begin{equation*}
	    \frac{q_{b_n}^{2\delta}}{q_{b_n}^{4}} \sum_{0 < |q| < q_{b_n}} \frac{1}{\Vert q \alpha \Vert^{2}} \ll q_{b_n}^{2\delta - 2} < q_{b_n}^{-\lambda(1+\delta)} \leq r_{b_n}^{-\lambda},
	\end{equation*}
	and this finishes the proof.
	
	\item[\textbf{\textit{Case 2:}}] \textit{For all sufficiently large $n$, we have $q_{n+1} < q_n^2$.}
	
	In this case we stick with the original sequence $\{r_n\}_{n \geq 1}$ and observe that for any $0 < k < n$ we can rewrite the sum in the RHS of \eqref{term_4_eq} as
	\begin{equation}\label{split_eq}
	\begin{split}
	    \left[\sum_{0 < |q| < q_k} + \sum_{q_k \leq |q| < q_n}\right] \frac{1}{|q|^{2\eps}} \frac{1}{\Vert q \alpha \Vert^{2}} &< \sum_{0 < |q| < q_k} \frac{1}{\Vert q\alpha \Vert^2} + q_k^{-2\eps} \sum_{q_k \leq |q| < q_n} \frac{1}{\Vert q\alpha \Vert^2} \\
	    & \ll q_k^2 + q_k^{-2\eps} q_n^2,
	\end{split}
	\end{equation}
	where once again we have used \autoref{cont_frac_whole_lemma}.
	
	Take $0<k<n$ such that $q_k \in \left[q_n^{1/4}, q_n^{1/2}\right]$, which exists for all $n$ sufficiently large since we can find such terms in any interval of the form $[a, a^2]$ for $a$ sufficiently large, because of the given condition. Then the corresponding upper bound when we plug \eqref{split_eq} into \eqref{term_4_eq} is
	\begin{equation*}
	    \frac{q_n^{2\delta}}{q_{n+1}^2} \left(q_n + q_n^{2-\eps/2}\right) \ll q_n^{2\delta - \eps/2} < q_n^{-\lambda(1+\delta)} \leq r_n^{-\lambda},
	\end{equation*}
	which establishes the result of \autoref{PR_lemma}.
\end{itemize}
    
\end{boldproof}


\section{Counterexample to polynomial rate of rigidity in \texorpdfstring{$C^1$}{C1}}\label{counterexample_section}

\autoref{PR_lemma} raises the question of how low one can push the smoothness of $\phi$ and still have a polynomial rate of rigidity for $T_{\alpha, \phi}$. We show that, at least along the sequence $\{q_n\}_{n \geq 1}$ of denominators of best rational approximations for an irrational $\alpha$, there is $\phi \in C^1$ with $\widehat{\phi}(0) = 0$ such that 
\begin{equation*}
    \int_{\T \times \T} d(T_{\alpha, \phi}^{q_n}(x, y), (x, y))^2 \dd\nu(x, y) \gg_\delta q_n^{-\delta}
\end{equation*}
for every $\delta > 0$, unlike what happens for $\phi \in C^{1+\eps}$ with $\widehat{\phi}(0)=0$ (observe that in that case $\ell_n = 1$ works in \autoref{PR_lemma}). 

Indeed, let $\alpha \in \R \setminus \Q$ and choose $\phi: \T \to \R$ given by
\begin{equation*}
	\phi(x) := \frac{1}{C} \sum_{k \geq 2} \frac{e(q_k x) + e(- q_k x)}{q_k (\log{q_k})^{2}},
\end{equation*}
where $C>0$ will be chosen to be sufficiently large. Since $q_k \geq 2^{(k-1)/2}$ by \ref{P1_item}, $\sum_{k \geq 2} (\log q_k)^{-2}$ is absolutely convergent and therefore $\phi \in C^1$. Take $C>0$ large enough so that $\Var(\phi) < 1/2$. By the Denjoy-Koksma inequality, we have 
\begin{equation*}
	\left|\sum_{j = 0}^{q_n - 1}\phi(x + j\alpha) - q_n \int_{\T} \phi(z) \dd z \ \right| \leq \Var(\phi) < \frac{1}{2}
\end{equation*}
for every $x \in \T$. Since $\widehat{\phi}(0)=0$ we conclude that $|S_{q_n}(\phi)(x)| < 1/2$ , so that $\Vert S_{q_n}(\phi)(x)\Vert = |S_{q_n}(\phi)(x)|$ for all $x \in \T$. Therefore, the beginning of the proof of \autoref{PR_lemma} shows that
\begin{equation}\label{counterexample_eq}
    \int_{\T \times\T} d(T_{\alpha, \phi}^{q_n}(x, y), (x, y))^2 \dd \nu(x, y) \asymp \Vert q_n \alpha \Vert^2 + \sum_{k \geq 2} \frac{1}{q_k^2 (\log{q_k})^4} \left| \frac{1 - e(q_n q_k \alpha)}{1 - e(q_k \alpha)}\right|^2.
\end{equation}

If $q_n \cdot \Vert q_n \alpha \Vert < 1/2$ then $\Vert q_n^2 \alpha \Vert = q_n \cdot \Vert q_n \alpha \Vert$, so
\begin{equation*}
    \frac{1}{q_n^2 (\log{q_n})^4}\left| \frac{1 - e(q_n^2 \alpha)}{1 - e(q_n \alpha)}\right|^2 \asymp \frac{1}{q_n^2 (\log{q_n})^4}\frac{\Vert q_n^2 \alpha \Vert^2}{\Vert q_n \alpha \Vert^2} = \frac{1}{(\log{q_n})^4}.
\end{equation*}

If instead $q_n \cdot \Vert q_n \alpha \Vert > 1/2$ then from $\Vert q_n \alpha \Vert < 1/q_{n+1}$ (by \ref{P2_item}) we get $q_{n+1} < 2 q_n$. Since $q_{n+2}> 2 q_n$ we have $q_n \cdot \Vert q_{n+1} \alpha \Vert < q_n/q_{n+2}< 1/2$, so $\Vert q_n q_{n+1} \alpha \Vert = q_n \cdot \Vert q_{n+1} \alpha \Vert$, and in conclusion
\begin{equation*}
    \frac{1}{q_{n+1}^2 (\log{q_{n+1}})^4}\left| \frac{1 - e(q_n q_{n+1} \alpha)}{1 - e(q_{n+1} \alpha)}\right|^2 \asymp \frac{1}{q_n^2 (\log{q_n})^4}\frac{\Vert q_n q_{n+1} \alpha \Vert^2}{\Vert q_{n+1} \alpha \Vert^2} = \frac{1}{(\log{q_n})^4}.
\end{equation*}

Taking respectively the terms corresponding to $k = n$ and $k = n+1$ in \eqref{counterexample_eq} and using positivity of the other terms we conclude that the whole expression is $\gg (\log{q_n})^{-4}$, so there is no polynomial rate of convergence to zero along any subsequence of $\{q_n\}_{n \geq 1}$. In fact, \cite[Lemma 3.2]{C2} shows that a decay of the form $\exp(-(\log{\log{q_n}})^{1+\delta})$ for any $\delta > 0$ would be enough for M\"obius disjointness, but that too is false by our counterexample.


\section{Extension of general rigidity results to \texorpdfstring{$\phi$}{phi} of non-zero mean}\label{abs_con_section}

Recall that a topological dynamical system $(X, T)$ is called \textit{rigid} if for each $\nu \in M(X, T)$ there exists a sequence $\{r_n\}_{n\geq 1}$ of positive integers such that $g \circ T^{r_n} \to g$ in $L^2(\nu)$ for all $g  \in L^2(\nu)$. 

By theorems of Herman \cite[XIII.4.8]{herman_abs_con} and Gabriel, Lema\'{n}czyk, and Liardet \cite[Th{\'e}or{\`e}me 1.1]{uniform_rigidity}, if $\alpha \in \R \setminus \Q$ and $\phi$ is absolutely continuous, has topological degree zero, and satisfies $\widehat{\phi}(0) = 0$, then the skew product $T_{\alpha, \phi}$ is rigid, and in fact they show that $T_{\alpha, \phi}^{q_n} \to \Id$ uniformly by obtaining 
\begin{equation*}
    \lim_{n \to \infty} \sup_{x \in \T} |S_{q_n}(\phi)(x)| = 0.
\end{equation*} 

Lema\'{n}czyk and Mauduit \cite[Theorem 1]{rigidity_L_2} (see also \cite[Corollary 2.8]{koksma_banach}) generalized\footnote{If $\phi: \T \to \R$ is absolutely continuous then $\phi' \in L^1(\T)$, so the Riemann-Lebesgue lemma gives $\widehat{\phi}(m) = o(1/|m|)$.} these theorems to show rigidity (though not uniformly) of $T_{\alpha, \phi}$ for all $\alpha \in \R \setminus \Q$ and $\phi \in L^2(\T)$ (of topological degree zero) satisfying $\widehat{\phi}(0) = 0$ and $\widehat{\phi}(m) = o(1/|m|)$.

The techniques of this paper may be employed to extend both results to cover the case $\widehat{\phi}(0) \not= 0$. Furthermore, in the case $\phi \in C^{1+\eps}$ we can modify \autoref{PR_lemma} to recover a uniform polynomial rate of rigidity instead of just the result in $L^2(\nu)$ presented previously.

\subsection{Uniform rigidity for \texorpdfstring{$\phi$}{phi} absolutely continuous}

\begin{proposition}\label{abs_con_prop}
If $\alpha \in \R \setminus \Q$ and $\phi$ is absolutely continuous of topological degree zero, then the skew product $T_{\alpha, \phi}$ is uniformly rigid.
\end{proposition}

\begin{boldproof}[Proof of \autoref{abs_con_prop}]
    We can simply use the original result for the zero mean case to conclude that there is $\lambda(n) \to \infty$ as $n \to \infty$ such that
\begin{equation*}
    \sup_{x \in \T} |S_{q_n}(\phi - \widehat{\phi}(0))(x)| \leq \lambda(n)^{-1},
\end{equation*}
so choose $\ell_n \in \Z$ with
\begin{equation*}
	0 < \ell_n \leq \lambda(n)^{1/2} \quad \text{and} \quad \Vert \ell_n q_n \widehat{\phi}(0) \Vert < \lambda(n)^{-1/2}
\end{equation*}
using Dirichlet's approximation theorem to get
\begin{equation*}
\begin{split}
    \Vert S_{\ell_n q_n}(\phi)(x) \Vert &\leq \Vert \ell_n q_n \widehat{\phi}(0) \Vert + |S_{\ell_n q_n}(\phi - \widehat{\phi}(0))(x)| \\
    &< \lambda(n)^{-1/2} + \sum_{k = 0}^{\ell_n - 1} |S_{q_n}(\phi - \widehat{\phi}(0))(x + k q_n \alpha)| \ll \lambda(n)^{-1/2} \to 0
\end{split}
\end{equation*}
uniformly in $x \in \T$. Therefore, $T_{\alpha, \phi}$ is uniformly rigid along the sequence $\{\ell_n q_n\}_{n\geq 1}$.

\end{boldproof}

\subsection{Rigidity for \texorpdfstring{$\phi$}{phi} with tamely decaying Fourier coefficients}

\begin{proposition}\label{tame_decay_prop}
If $\alpha \in \R \setminus \Q$ and $\phi \in L^1(\T)$ (of topological degree zero) satisfies $\widehat{\phi}(m) = o(1/|m|)$, then the skew product $T_{\alpha, \phi}$ is rigid.
\end{proposition}

\begin{boldproof}[Proof of \autoref{tame_decay_prop} (Sketch)]
    We have the Fourier expansion\footnote{It follows from the conditions that $\phi \in L^2(\T)$.} (in $L^2(\T)$)
\begin{equation*}
	\phi(x) = \sum_{q \in \Z} c_q e(qx) \quad \text{with}\quad |c_q| \leq \frac{1}{|q| \cdot \psi(|q|)} \ \text{for} \  q \not= 0,
\end{equation*}
where $\psi:(0, \infty) \to (0, \infty)$ satisfies $\psi(z) \to \infty$ as $z \to \infty$ and for technical reasons we can of course also assume that it is non-decreasing and does not grow too fast, say $\psi(z) \ll_\phi z^{1/100}$.

With the conditions above, we can show that there is a sequence of positive integers $\{r_n\}_{n\geq 1}$ such that
\begin{equation*}
	\int_{\T \times \T} d(T_{\alpha, \phi}^{r_n}(x, y), (x, y))^2 \dd \nu(x, y) \ll_\phi \psi(q_n^{1/4})^{-1/100} \to 0 \quad \text{as} \quad n \to \infty
\end{equation*}
for any $\nu \in M(\T^2, T_{\alpha, \phi})$.

The proof is a simple modification of the proof of \autoref{PR_lemma}, substituting $q_n^\eps$ with $\psi(q_n)$, so for instance $\ell_n \in \Z$ is chosen so that
\begin{equation*}
	0 < \ell_n \leq \psi(q_n)^{1/10} \quad \text{and} \quad \Vert \ell_n q_n c_0 \Vert < \psi(q_n)^{-1/10}.
\end{equation*}
Observe that we do not have multiplicativity of $\psi$, which is why the bound is not of the form $\psi(r_n)^{-1/100}$, but it is enough to prove that $T_{\alpha, \phi}$ is rigid\footnote{The bound implies rigidity for $T_{\alpha, \phi}$ since the Lipschitz continuous functions on $\T^2$ are dense in $L^2(\nu)$, for any $\nu \in M(\T^2, T_{\alpha, \phi})$. This follows by the Stone-Weierstrass theorem and the fact that $C(\T^2)$ is dense in $L^2(\nu)$, since $\nu$ is a Radon measure -- see for instance \cite[Proposition 7.9]{folland}.} (the latter bound could be obtained if we imposed extra attainable conditions on $\psi$).

\end{boldproof}

\subsection{Uniform polynomial rate of rigidity for \texorpdfstring{$\phi \in C^{1+\eps}$}{phi in C(1 + epsilon)}}\label{uniform_poly_rigidity_subsec}

Finally, we point out that the conclusion of \autoref{PR_lemma} can actually be strengthened to a uniform polynomial rate of rigidity:

\begin{proposition}\label{uniform_rigidity_prop}
 Let $0 < \eps < \frac{1}{100}$ and $\alpha \in \R \setminus \Q$. If $\phi: \T \to \R$ is of class $C^{1 + \eps}$, then there exists an unbounded sequence of positive integers $\{r_n\}_{n\geq1}$ such that
\begin{equation*}
    \sup_{(x, y) \in \T \times \T} d(T_{\alpha, \phi}^{r_n}(x, y), (x, y)) \ll_{\phi, \eps} r_n^{-\eps/200}.
\end{equation*}
\end{proposition}

\begin{boldproof}[Proof of \autoref{uniform_rigidity_prop} (Sketch)]
We start by substantially modifying the results of \autoref{cont_frac_whole_lemma} and \autoref{cont_frac_slice_lemma}. Namely one can show, using the same techniques as in the corresponding results of \autoref{cont_frac_section} but this time for the functions
\begin{equation*}
        f_1(z) = 
        \begin{cases}
        2 q_k, & \text{if} \ \Vert z \Vert \leq \frac{1}{2q_k} \\
        \frac{1}{\Vert z \Vert}, & \text{if} \ \Vert z \Vert > \frac{1}{2q_k}
        \end{cases}
        \quad \text{and} \quad
        f_2(z) = 
        \begin{cases}
        c, & \text{if} \ \Vert z \Vert \leq \frac{1}{c} \\
        \frac{1}{\Vert z \Vert}, & \text{if} \ \Vert z \Vert > \frac{1}{c},
        \end{cases}
    \end{equation*}
respectively, that if $\eps > 0$, $\alpha \in \R \setminus \Q$, $k \geq 1$, and $1 \leq c \leq q_k$ then
\begin{equation}\label{cont_frac_whole_2_eq}
	\sum_{0 < |q| < q_k} \frac{1}{\Vert q \alpha \Vert} \ll q_k \log(q_k+1)
\end{equation}
and 
\begin{equation}\label{cont_frac_slice_2_eq}
	\sum_{q_k \leq |q| < q_{k+1}} \frac{1}{|q|^{1+ \eps}} \min \left\{\frac{1}{\Vert q \alpha \Vert}, c \right\} \ll_\eps  \frac{\log(c+1)}{q_k^\eps}.
\end{equation}

Then expanding $S_{r_n}(\phi)(x)$ into a Fourier series and trivially bounding it we get
\begin{equation*}
    d(T_{\alpha, \phi}^{r_n}(x, y), (x, y)) \leq \Vert r_n \alpha \Vert + \Vert c_0 r_n \Vert + \sum_{q\not= 0} |c_q| \left| \frac{1 - e(q r_n \alpha)}{1 - e(q \alpha)} \right|,
\end{equation*}
so we can proceed as in the proof of \autoref{PR_lemma} with the expression above corresponding to \eqref{square_eq} and the bounds of \eqref{cont_frac_whole_2_eq} and \eqref{cont_frac_slice_2_eq} corresponding to \autoref{cont_frac_whole_lemma} and \autoref{cont_frac_slice_lemma}, respectively, to get the desired uniform polynomial decay. 
    
\end{boldproof}

\begin{remark}
The proof actually shows that for every $\phi: \T \to \R$ of class $C^{1+\eps}$ and of mean zero,
\begin{equation}\label{uniform_poly_S_bound_eq}
    \sup_{x \in \T} | S_{q_n}(\phi)(x)| \ll_{\phi, \eps} q_n^{-\eps/200},
\end{equation}
since in that case we can take $\ell_n = 1$ throughout the argument.
\end{remark}

\begin{remark}
Even though \autoref{uniform_rigidity_prop} gives a stronger result than \autoref{PR_lemma}, we chose to emphasize the latter in our presentation because the $L^2$ methods employed there seem more suitable for generalization (and the proof is slightly more complicated). For instance, an approach to \autoref{tame_decay_prop} using $L^\infty$ methods would already be frustrated by the presence of the extra logarithmic factor in \eqref{cont_frac_whole_2_eq}, if the decay of the Fourier coefficients is sufficiently slow. Therefore, the use of $L^2$ methods seems to allow us to go a bit further.
\end{remark}


\section{Smooth flows on \texorpdfstring{$\T^2$}{T2} and Rokhlin extensions}\label{flows_rokhlin_section}

We can adapt the result of this paper, following \cite{C2}, to give M{\"o}bius disjointness for new cases of smooth flows on the torus and Rokhlin extensions.

\subsection{Smooth flows on \texorpdfstring{$\T^2$}{T2}}

For $\alpha \in \R\setminus \Q$, let $f :\T \to \R$ be a strictly positive continuous function. Let 
\begin{equation*}
    \T^f := \{(x, s) \in \T \times \R : 0 \leq s \leq f(x)\}/\sim
\end{equation*}
where $\sim$ denotes the equivalence relation $(x, s + f(x)) \sim (R_\alpha(x), s)$ in $\T \times \R$ and $R_\alpha : \T \to \T$ is the irrational rotation by $\alpha$. We can define a \textit{special flow} $T^f = \{T^{f}_t\}_{t\in \R}$ over $R_\alpha$ with roof function $f$, which acts on $\T^f$ by
\begin{equation*}
    T^f_t(x, s) := (x, s + t)
\end{equation*}
for all $(x, s)\in \T^f$. More explicitly, if we extend a previous definition to
\begin{equation*}
    S_N(f)(x) := 
    \begin{cases}
        \sum_{0 \leq j < N} f(R_\alpha^j(x)), & \text{if} \ N > 0 \\
        0, & \text{if} \ N = 0 \\
        \sum_{N \leq j < 0} f(R_\alpha^j(x)), & \text{if} \ N < 0
    \end{cases}
\end{equation*}
then
\begin{equation*}
    T^f_t(x, s) = (R_\alpha^N(x), s + t - S_{N}(f)(x))
\end{equation*}
for all $(x, s)\in \T^f$, where $N = N(x, s, t) \in \Z$ is such that 
\begin{equation*}
    S_{N}(f)(x) \leq s+t < S_{N+1}(f)(x),
\end{equation*} 
which exists and is unique as $f$ is continuous and strictly positive.

Every sufficiently smooth area-preserving flow on $\T^2$ with no fixed points or closed orbits can be represented by such a special flow for $f$ with corresponding smoothness properties (see \cite{ergodic_book}).

We have the following consequence of our work:

\begin{corollary}\label{smooth_flow_cor}
Let $\eps > 0$ and $\alpha \in \R\setminus\Q$. If $f \in C^{1+\eps}$ then all the maps of the special flow $T^f = \{T_t^f\}_{t\in\R}$ over the irrational rotation $R_\alpha$ are M{\"o}bius disjoint.
\end{corollary}

\begin{boldproof}[Proof of \autoref{smooth_flow_cor}]
    There is a natural quotient metric $D$ making $\T^f$ a compact metric space (see \cite[Appendix 9.1]{quotient_metric}), and it satisfies
\begin{equation}\label{metric_ineq_eq}
    D(T^f_t(x, s), (x, s)) \leq |t| \quad \text{for all} \quad t\in \R \quad \text{and} \quad (x, s) \in \T^f.
\end{equation}

Denote $\beta := \widehat{f}(0) > 0$ and let $q_n$ be the denominators of convergents of the continued fraction of $\alpha$, as before. For a fixed $t \in \R$, let $v_n \in \Z$ be such that
\begin{equation*}
    0 < v_n \leq q_n^{1 + \gamma} \quad \text{and} \quad \left\Vert v_n \frac{t}{q_n \beta} \ \right\Vert < q_n^{-1-\gamma},
\end{equation*}
where $\gamma >0$ will be chosen later to be sufficiently small ($v_n$ exists by Dirichlet's approximation theorem). Then there is $j_n \in \Z$ such that 
\begin{equation}\label{tv_n_eq}
    |t v_n - j_n q_n \beta| < \frac{q_n \beta}{q_n^{1+\gamma}} \ll_f q_n^{-\gamma}
\end{equation}
and
\begin{equation}\label{j_n_mod_eq}
    |j_n| < \left|\frac{v_nt}{q_n \beta} \right| + 1 \ll_{f, t} q_n^\gamma.
\end{equation}

For every $(x, s) \in \T^f$ we have
\begin{equation*}
\begin{split}
    D(T^f_{t v_n}(x, s), (x, s)) &= D(T^f_{t v_n - S_{j_n q_n}(f)(x)}\circ T^f_{S_{j_n q_n}(f)(x)}(x, s), (x, s)) \\
    &\leq D(T^f_{t v_n - S_{j_n q_n}(f)(x)}(R_\alpha^{j_n q_n}(x), s), (R_\alpha^{j_n q_n}(x), s)) + D((R_\alpha^{j_n q_n}(x), s), (x, s)) \\
    &\leq |S_{j_n q_n}(f - \beta)(x)| + |tv_n - j_n q_n \beta| + \Vert j_n q_n \alpha \Vert \\
    & \ll_{f, t} q_n^\gamma \cdot \sup_{z\in \T} |S_{q_n}(f-\beta)(z)| + q_n^{-\gamma} + q_n^\gamma \cdot \Vert q_n \alpha \Vert,
\end{split}
\end{equation*}
where we have used \eqref{metric_ineq_eq}, \eqref{tv_n_eq} and \eqref{j_n_mod_eq}. Choosing $\gamma := \eps/1000$ and using \eqref{uniform_poly_S_bound_eq} (we could also take the $L^2$ norm and use the proof of \autoref{PR_lemma}) we get the bound $\ll_{f, t, \eps} q_n^{-\eps/1000} < v_n^{-\eps/2000}$, which gives a polynomial rate of rigidity for $T_t^f:\T^f \to \T^f$ along the (unbounded, unless $t = 0$) sequence $\{v_n\}_{n\geq 1}$, and this implies M{\"o}bius disjointness for $T_t^f$.

\end{boldproof}

\subsection{Rokhlin extensions} As before, let $\alpha \in \R \setminus \Q$ and let $R_\alpha: \T \to \T$ denote the irrational rotation by $\alpha$. Given a continuous function $f: \T \to \R$, a compact metric space $(Y, \rho)$ and a continuous flow $L = \{L_t\}_{t\in\R}$ acting on $Y$, we can define a \textit{Rokhlin extension} $E_{f, L}$ of $R_\alpha$, acting on $\T \times Y$ by
\begin{equation*}
    E_{f, L}(x, y) := (R_\alpha(x), L_{f(x)}(y))
\end{equation*}
for all $(x, y) \in \T \times Y$ (observe that if $Y = \T$ and $L$ is the linear flow we recover the Anzai skew product $T_{\alpha, f}$). We have the following disjointness result in this case:

\begin{corollary}\label{rokhlin_con}
Let $\eps > 0$. If $f \in C^{1+\eps}$ has mean zero and the flow $L = \{L_t\}_{t\in \R}$ is uniformly Lipschitz continuous in $t$, then $E_{f, L}$ is M{\"o}bius disjoint.
\end{corollary}

\begin{boldproof}[Proof of \autoref{rokhlin_con}]
    If $D$ denotes the product metric in $\T \times Y$ then
\begin{equation*}
    D(E_{f, L}^{q_n}(x, y), (x, y)) = \Vert q_n \alpha \Vert + \rho(L_{S_{q_n}(f)(x)}(y), y) \ll_L \Vert q_n \alpha \Vert + |S_{q_n}(f)(x)|,
\end{equation*}
where the implied constant does not depend on $(x, y)$. Using \eqref{uniform_poly_S_bound_eq} we get a polynomial rate of rigidity for $E_{f, L}$ along $\{q_n\}_{n\geq 1}$ (we could also take the $L^2$ norm and use the proof of \autoref{PR_lemma}), so the corollary follows. 

\end{boldproof}


\nocite{*}  
\bibliographystyle{abbrv}
\bibliography{references}

\end{document}